\documentclass[12pt,reqno]{amsart}

\usepackage{amsmath, amscd, amssymb, amsthm,amsfonts}
\usepackage{euscript, mathrsfs, latexsym,mathabx,stmaryrd}
\usepackage{color}

\usepackage{lmodern}
\usepackage[T1]{fontenc} 

\usepackage{multicol}

\ifx\pdfoutput\undefined
\usepackage{graphicx}
\else
\fi

\ifx\xetex
\usepackage{fontspec}
\usefont{EU1}{lmr}{m}{n}
\else
\fi

\usepackage{url} 
\usepackage{tikz}
\usetikzlibrary{matrix,arrows,positioning}
\tikzset{node distance=2cm, auto}

\usepackage[nodayofweek]{datetime}

\parskip=\medskipamount

 \DeclareMathOperator{\cofib}{cofib}
\DeclareMathOperator{\colim}{colim}

\newcommand{\conj}[1]{\quad\textnormal{ #1 }\quad}

\newcommand{\tl}[1]{\ensuremath{\tilde{#1}}}
\newcommand{\rd}[1]{\ensuremath{\{#1\}}}
\newcommand{\inp}[1]{\ensuremath{\langle #1 \rangle}}

\newcommand{\normaltext}[1]{\textnormal{#1}}

\makeatletter
\def\imod#1{\allowbreak\mkern2.5mu({\operator@font mod}\,#1)}
\makeatother

\newcommand{\ott}{\otimes}

\renewcommand{\d}{\delta}

\newcommand{\MPic}[1]{
\begin{minipage}{.35in}
\includegraphics[scale=.45]{#1}
\end{minipage}
}

\newcommand{\aX}{\mathcal{X}}
\newcommand{\aY}{\mathcal{Y}}

\newcommand{\bA}{\mathbf{A}}

\usepackage{bbm}

\newcommand{\FF}{\mathbb{F}}

\newcommand{\ZZ}{\mathbb{Z}}

\newcommand{\kk}{\mathbb{k}}

\theoremstyle{plain}

\newtheorem{thm}{Theorem}[section]

\newtheorem{theorem}[thm]{Theorem}

\newtheorem{prop}[thm]{Proposition}

\newtheorem{corollary}[thm]{Corollary}

\newtheorem{lemma}[thm]{Lemma}

\theoremstyle{remark}

\theoremstyle{definition}

\newtheorem{definition}[thm]{Definition}

\newtheorem{remark}[thm]{Remark}

\numberwithin{equation}{section}

\usetikzlibrary{knots}
\usepackage[colorlinks,linkcolor=blue,citecolor=blue]{hyperref}
\usepackage{overpic}
\usepackage{enumitem}   

\begin{document}

\title{An incompatibility between spectrification and the Szab\'{o} spectral sequence}
\date{\today}
\author[Benjamin Cooper]{Benjamin Cooper}
\author[Pravakar Paul]{Pravakar Paul}
\author[Nicholas Seguin]{Nicholas Seguin}
\address{University of Iowa, Department of Mathematics, 14 MacLean Hall, Iowa City, IA 52242-1419 USA}
\email{ben-cooper\char 64 uiowa.edu}
\email{pravakar-paul\char 64 uiowa.edu}
\email{nicholas-seguin\char 64 uiowa.edu}

\def\JS#1{\textcolor[rgb]{0,.75,.8}{ [JS: #1]}}
\def\BC#1{\textcolor[rgb]{0,.5,1}{ [BC: #1]}}
\def\OLD#1{\textcolor[rgb]{1,.69,.4}{ [OLD STUFF: #1]}}
\def\OPT#1{\textcolor[rgb]{1,.61,0}{ [OPTIONAL: #1]}}

\newcommand{\tm}{\widetilde{m}}
\newcommand{\dga}[0]{\operatorname{-dga}}
\newcommand{\dgca}[0]{\operatorname{-dgcoa}}
\newcommand{\Om}{\Omega}
\newcommand{\pa}{\partial}

\newcommand{\Ainf}{A_\infty}
\newcommand{\varB}[1]{{\operatorname{\mathit{#1}}}}
\renewcommand{\slash}[1]{H_/(#1)}
\newcommand{\pAinf}[0]{\varB{p-\Ainf}\!}
\newcommand{\semp}{{\{\emptyset\}}}
\newcommand{\z}{z}
\newcommand{\T}{T}
\renewcommand{\d}{\delta}
\newcommand{\h}{h}
\renewcommand{\t}{t}
\newcommand{\G}{\Gamma} 
\newcommand{\A}{\Lambda} 
\renewcommand{\bA}{\bar{\Lambda}} 

\renewcommand{\kk}{k}

\newcommand{\lab}[1]{\normaltext{#1}}
\newcommand{\nt}[1]{\normaltext{#1}}
\newcommand{\Ya}{\mathcal{Y}\hspace{-.1475em}a_{2,2}}
\newcommand{\vnp}[1]{\lvert #1 \rvert}
\newcommand{\vvnp}[1]{\lvert\lvert #1 \rvert\rvert}
\newcommand{\I}{[0,1]}
\newcommand{\xto}[1]{\xrightarrow{#1}}
\newcommand{\xfrom}[1]{\xleftarrow{#1}}
\newcommand{\from}{\leftarrow}

\newcommand{\dgcat}{\normaltext{dgcat}_k}
\newcommand{\Forget}{\normaltext{Forget}}
\newcommand{\Mat}{Mat}
\newcommand{\op}{\normaltext{op}}
\newcommand{\ab}{\normaltext{ab}}

\newcommand{\Vect}{Vect}
\newcommand{\Kom}{Kom}
\newcommand{\Set}{Set}
\newcommand{\Ch}{Ch}
\newcommand{\Hom}{Hom}
\renewcommand{\d}{\delta}

\newcommand{\pre}[0]{\operatorname{pre-}}
\newcommand{\coker}[0]{\operatorname{coker}}
\newcommand{\im}[0]{\operatorname{im}}


\begin{abstract}
The Lipshitz-Sarkar Steenrod operations on (even) Khovanov homology are
incompatible with Szab\'{o}'s geometric spectral sequence. Obstructions to
integral lifting and spectrification are observed.
  \end{abstract}

\maketitle

\section{Introduction}
Variations on constructions of knot homology theories have determined a rich
mathematical structure, reflecting both their many origins and the
surprising capacity for dissimilar settings to carry the complexity of knots
and links. In the presence of many different, but ultimately equivalent
definitions, it is common to use the extra structure available in one
setting to predict new structure in another.  Here we show that extra
structure present in one setting is {\em incompatible} with extra structure
from another.

R. Lipshitz and S. Sarkar introduced a spectrum-level refinement of Khovanov
homology which determines Steenrod operations $Sq^n :
Kh^{t,q}(K;\FF_2)\to Kh^{t+n,q}(K;\FF_2)$ for $n\geq 1$ \cite{LS, LSS}.
Z. Szab\'{o} defined a spectral sequence which interpolates from Khovanov
homology to a knot homology theory (conjecturally) isomorphic to the
Heegaard-Floer homology $\widehat{HF}(-\Sigma(K)\#(S^1\times S^2))$ of the
double branched cover $\Sigma(K)$ connected sum $S^1\times S^2$ \cite{OZ,
  Szabo}.  We prove that for each $n\geq 1$, there is a knot $K_n$ so that
$Sq^n$ does not commute with the differentials of Szab\'{o}'s spectral
sequence.

Since compatibility between these operations and the spectral sequence is
implied by several constructions: a spectrification of Szab\'{o} homology
admitting a filtration from which the (even) Khovanov spectrum can be
recovered or an integral lift of the Szab\'{o} spectral sequence from the
(even) Khovanov homology, our counterexample shows that these constructions
are not possible.


\subsection*{Notation}  
Unless indicated otherwise, coefficients are given by the field with two elements $\FF_2$ and chain complexes $(C,d)$ will be cochain complexes; $\vnp{d}_t = +1$.

\section{Spectral sequences}
Suppose that $(CKh(K),d)$ is the Khovanov chain complex associated to a link
$K$\!. A map $\d : CKh(K) \to CKh(K)$ of chain complexes is a {\em
  perturbation} when $(d+\d)^2=0$.  When the bidegree of $\d$ differs from
the bidegree of $d$, there is a filtration of $(CKh(K),d+\d)$ with
associated graded given by the original complex $(CKh(K),d)$. And so,
associated to a perturbation $\d$, there is a spectral sequence $\{E_n,d_n\}_{n=2}^{\infty}$, consisting of pages $E_n = \oplus_{(t,q)\in\ZZ\times\ZZ} E^{t,q}_n$ and differentials $d_n : E_n \to E_n$, such that $E_{n+1} := H(E_n, d_n)$, from the Khovanov homology $E_2 = Kh(K)$, (where $Kh(K) := H(CKh(K), d)$), to the homology of the total complex $E_\infty = H(CKh(K),d+\d)$ \cite[Thm. 2.6]{McCleary}.

\begin{definition}\label{risetoopdef}
An endomorphism $f_m : E_m \to E_m$ acting on the $E_m$-page {\em gives rise to an operation on} $\{E_n\}_{n=m}^{\infty}$ when there is a sequence of linear maps  $f_n : E_n \to E_n$  for  $n>m$ which satisfy the two properties below.
\begin{enumerate}
\item The map $f_n : E_n \to E_n$ commutes with the differential on the $E_n$-page: $f_nd_n = d_n f_n$.
\item The map $f_{n+1} : E_{n+1} \to E_{n+1}$ agrees with the map induced by $f_n$ on the $E_n$-page: $f_{n+1} = H(f_n, d_n) : H(E_n, d_n) \to H(E_n, d_n)$.
\end{enumerate}
\end{definition}

Recall that the Khovanov chain complex $(CKh(K), d)$ is defined as a direct sum of
tensor products of the Frobenius algebra $\FF_2[x]/(x^2)$. Choosing any
point on the knot $K$ allows us to adjoin a handle corresponding to
multiplication by $x$. This observation gives rise to a chain map $X : CKh(K) \to CKh(K)$ and an induced map $X_*$ on homology \cite[\S 3]{KPatterns}.
A. Shumakovitch introduced a decomposition of the form
\begin{equation}\label{reducedeq}
Kh(K;\FF_2) \cong \widetilde{Kh}(K;\FF_2) \oplus X_*\widetilde{Kh}(K;\FF_2)
\end{equation}
where $\widetilde{Kh}(K;\FF_2)\subset Kh(K;\FF_2)$ and $\widetilde{Kh}(K;\FF_2) \cong Kh(K;\FF_2)\ott_{\FF_2[x]/(x^2)} \FF_2$ is the reduced Khovanov homology \cite[\S 3]{Shum}.

The lemma below shows that the first part of this story extends to Szab\'{o}'s spectral sequence.

\begin{lemma}
The map $X_* : Kh(K) \to Kh(K)$ gives rise to an operation on the Szab\'{o}'s spectral sequence.
\end{lemma}
\begin{proof}
We construct a map $X_{Sz} : CKh(K)\to CKh(K)$ which commutes with Szab\'{o}'s differential $\d_{Sz}$ and agrees with $X$ on the associated graded of the filtration defining the spectral sequence. This is accomplished by refactoring Szab\'{o}'s proof of invariance under the Reidemeister 1 move \cite[Thm. 7.2]{Szabo}.

Pick a point $p$ on $K$. Adding a kink in the knot at $p$ gives a knot diagram $K^\sharp$. By resolving the crossing at the kink, the chain complex $\inp{K^\sharp}$
\begin{center}
\begin{tikzpicture}[scale=10, node distance=2.5cm]
\node (X) {$\left[\MPic{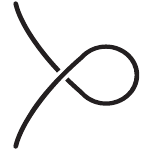}\,\,\,\right]_{(\inp{K^\sharp},\d_{Sz})}$};
\node (Y) [right=.3cm of X] {$=$};
\node (Z) [right=.3cm of Y] {$Cone\Bigg(\left[\MPic{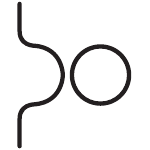}\,\,\right]_{(C_0,\d_0)}$};
\node (A) [right=1.5cm of Z] {$\left[\MPic{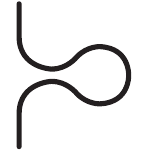}\,\,\right]_{(C_1,\d_1)}\Bigg)$};
\draw[->] (Z) to node {$S$} (A);
\end{tikzpicture} 
\end{center}
associated to the diagram $K^{\sharp}$, can be written as a cone on a map $S$, so $\d_{Sz}^2=0$ implies $S\d_0 = \d_1 S$. Now the disjoint circle in the diagram for $C_0$ produces a decomposition $C_0 \cong C_0^+ \oplus C_0^-$ where $C_0^-$ consists of elements divisible by $x$ (in the Frobenius algebra associated to the disjoint circle) and $C_0^+$ those elements which are not divisible by $x$. Under this isomorphism $S$ is a sum of two maps $S = X_{Sz}+1$ where $X_{Sz} : C_0^- \to C_1$ and $1 : C_0^+ \to C_1$. Szab\'{o} observes that the map $1$ is the identity map by construction. So $S$ and $1$ are chain maps, which implies that $X_{Sz}$ must also commute with Szab\'{o} differentials.

Since the map $X_{Sz}$ has positive homological degree or $t$-degree, it preserves the filtration $F^p C := \{ c\in C : \vnp{c}_t \geq p \}$ defining the Szab\'{o} spectral sequence and induces an operation on the spectral sequence \cite[Thm. 3.5]{McCleary}. Lastly, since the first order term of the Szab\'{o} differential is the Khovanov differential, the first order term of the map $X_{Sz}$ is $X$, so $X$ (and $X_* = H(X,d)$) extend to operations on the spectral sequence.

  \end{proof}

The map $X_{Sz}$ can be written as $X_{Sz} = \sum_{n=1}^{\infty} X_n$ where
$X_1 = X$ and $X_n := \sum_{p+q=n} E_{p,q}$ for $n\geq 1$ where $E_{p,q}$ is the
assignment from \cite[Def. 4.5]{Szabo}. This formula depends on an 
orientation of the saddle $S$, but any two choices are homotopic \cite[\S 5.1]{Szabo}.

The map $P(t)$ from the proof of Prop. 1.3 \cite{Seed} is also an extension of the map $X$. 


\section{The counterexample}
We now combine the materials above with the output of the computer programs
KnotKit by C. Seed and JavaKh by D. Bar-Natan and J. Green
\cite{Seed, SeedSq, fastkh} to produce an incompatibility between
the Bockstein $Sq^1 : Kh(K) \to Kh(K)$ and the Szab\'{o} spectral sequence. This occurs on the
$E_3$-page of the spectral sequence associated to the torus knot $T(4, 5)$.

The Poincar\'{e} polynomial of the unreduced $\FF_2$-Khovanov homology of the torus knot $T(4,5)$ is given by
\begin{align*}
  P_2 &= (q^{11} + q^{13}) t^0 + (q^{15}+q^{17})t^2 + (q^{17} + q^{19})t^3 + (q^{17} + q^{19})t^4 + (q^{21}+q^{23})t^5 \\
  & \quad\quad + (q^{19} + 2 q^{21} + q^{23})t^6 + (q^{21}+ 2q^{23} + q^{25}) t^7 + (q^{23} + q^{25}) t^8 \\
  & \quad\quad + (q^{25}+2 q^{27} + q^{29})t^9 + (q^{27}+q^{29}) t^{10}
\end{align*}
This is the Poincar\'{e} polynomial of the $E_2$-page of the Szab\'{o} spectral
sequence. The polynomials associated to the $E_3$ and $E_4 = E_{\infty}$
pages are given below.
\begin{align*}
  P_3 &= (q^{11} + q^{13}) t^0 + (q^{17} + q^{19})t^3 + (q^{19} + 2 q^{21} + q^{23})t^6 + (q^{21}+ q^{23}) t^7\\
  & \quad\quad + (q^{23} + q^{25}) t^8 + (q^{25}+2 q^{27} + q^{29})t^9 + (q^{27}+q^{29}) t^{10}\\
  P_4 &= (q^{11} + q^{13}) t^0 + (q^{19} +  q^{21})t^6 + (q^{21}+ q^{23}) t^7\\
  & \quad\quad + (q^{23} + q^{25}) t^8 + (q^{25}+2 q^{27} + q^{29})t^9 + (q^{27}+q^{29}) t^{10}
\end{align*}

The diagram below also contains this information. In the diagram, the non-zero $d_2$ differentials are denoted by solid arrows and non-zero $d_3$ differentials are denoted by dashed arrows.  The $(t,q)$-bidegree of the differential $d_n : E_n \to E_n$ is $(n,2n-2)$.

 \newcommand*\ecirc[1][1ex]{%
   \begin{tikzpicture}
   \draw (0,0) circle (.1);
   \end{tikzpicture}}
 \newcommand*\xcirc[1][1ex]{%
   \begin{tikzpicture}
   \fill (0,0) circle (.1);
   \end{tikzpicture}}

\begin{center}
\begin{tikzpicture}[scale=1]
  \draw[->] (0,0) -- (11.5,0) node[right] {$t$};
  \draw[->] (0,0) -- (0,10.5) node[above] {$q$};
  \draw[step=1] (0,0) grid (11,10);
  \draw (0.5,-.2) node[below] {$0$};
  \draw (1.5,-.2) node[below] {$1$};
  \draw (2.5,-.2) node[below] {$2$};
  \draw (3.5,-.2) node[below] {$3$};
  \draw (4.5,-.2) node[below] {$4$};
  \draw (5.5,-.2) node[below] {$5$};
  \draw (6.5,-.2) node[below] {$6$};
  \draw (7.5,-.2) node[below] {$7$};
  \draw (8.5,-.2) node[below] {$8$};
  \draw (9.5,-.2) node[below] {$9$};
  \draw (10.5,-.2) node[below] {$10$};
  \draw (-.2,0.5) node[left] {$11$};
  \draw (-.2,1.5) node[left] {$13$};
  \draw (-.2,2.5) node[left] {$15$};
  \draw (-.2,3.5) node[left] {$17$};
  \draw (-.2,4.5) node[left] {$19$};
  \draw (-.2,5.5) node[left] {$21$};
  \draw (-.2,6.5) node[left] {$23$};
  \draw (-.2,7.5) node[left] {$25$};
  \draw (-.2,8.5) node[left] {$27$};
  \draw (-.2,9.5) node[left] {$29$};

\fill (0.5, 0.5) circle (.15);
  \draw (0.5, 1.5) circle (.15);

  \fill (2.5, 2.5) circle (.15);
  \draw (2.5, 3.5) circle (.15);
\draw[->,double] (2.750, 3.5) -- (3.25, 3.5);
\draw[->,double, bend right=50] (2.750, 3.5) to (4.25, 3.5);

  \fill (3.5, 3.5) circle (.15);
  \draw (3.5, 4.5) circle (.15);

  \fill (4.5, 3.5) circle (.15);
  \draw (4.5, 4.5) circle (.15);

  \fill (5.5, 5.5) circle (.15);
\draw[->,double, bend right=50] (5.750, 5.5) to (7.25, 5.5);
  \draw (5.5, 6.5) circle (.15);
\draw[->,double, bend right=50] (5.750, 6.5) to (7.25, 6.5);
\draw[->,double, bend left=50] (5.750, 6.5) to (7.25, 6.5);

  \fill (6.5, 4.5) circle (.15);
  \draw (6.5, 5.5) node {$2$};
\draw[->, double] [bend left=20] (6.750, 5.5) to (7.25, 5.5);
\draw[->, double] [bend right=20] (6.750, 5.5) to (7.25, 5.5);
  \draw (6.5, 6.5) circle (.15);
\draw[->, double] [bend left=20] (6.750, 6.5) to (7.25, 6.5);
\draw[->, double] [bend right=20] (6.750, 6.5) to (7.25, 6.5);

  \fill (7.5, 5.5) circle (.15);
  \draw (7.5, 6.5) node {$2$};
  \draw (7.5, 7.5) circle (.15);

  \fill (8.5, 6.5) circle (.15);
  \draw (8.5, 7.5) circle (.15);

  \fill (9.5, 7.5) circle (.15);
  \draw (9.5, 8.5) node {$2$};
\draw[->, double] [bend left=20] (9.750, 8.5) to (10.25, 8.5);
\draw[->, double] [bend right=20] (9.750, 8.5) to (10.25, 8.5);
  \draw (9.5, 9.5) circle (.15);
\draw[->, double] (9.750, 9.5) -- (10.25, 9.5);

  \fill (10.5, 8.5) circle (.15);
  \draw (10.5, 9.5) circle (.15);

  \draw[->] (2.768, 2.634) -- (4.232, 3.366);
  \draw[->] (2.768, 3.634) -- (4.232, 4.366);
  \draw[->] (5.768, 5.634) -- (7.232, 6.366);
  \draw[->] (5.768, 6.634) -- (7.232, 7.366);

  \draw[->,dashed] (3.750, 3.666) -- (6.250, 5.334);
  \draw[->,dashed] (3.750, 4.666) -- (6.250, 6.334);
\end{tikzpicture}
\end{center}

The non-zero maps $Sq^1$ and $Sq^2$ of $(t,q)$-degrees $(1,0)$ and $(2,0)$
are represented by double arrows. A black dot is a generator which is in the
image of the operation $X_*$ and a white dot is a generator which is not in
the image of $X_*$, as in Eqn. \eqref{reducedeq}. The vector spaces of rank
2 are labelled by the number $2$. In the picture above, each such vector
space consists of one black dot and one white dot. Two horizontal arrows from a box labelled by the number $2$ are drawn to indicate that Steenrod operations are non-zero on both dots.

\begin{theorem}
  The Bockstein $Sq^1 : Kh(T(4,5)) \to Kh(T(4,5))$ does not give rise to an operation on the Szab\'{o} spectral sequence.
  \end{theorem}
\begin{proof}
  Assume that the Bockstein gives rise to an operation on the Szab\'{o} spectral sequence. Then there are maps
  $Sq^1_n : E^{t,q}_n  \to E^{t+1,q}_n$  for  $n\geq 2$
  which satisfy the two properties in the definition above. But this cannot be true! 

  First observe that, of the vector spaces: $E_2^{4,17}$, $E_2^{3,17}$, $E_2^{6,21}$ and $E_2^{7,21}$, only $E_2^{4,17}$ interacts with a non-zero $d_2$ differential. In this way, $E_3^{t,q} = H(E_2, d_2)$ implies the following isomorphisms
  $$E_3^{4,17} \cong 0 \conj{ and } E_3^{3,17} \cong E_2^{3,17},\quad E_3^{6,21} \cong E_2^{6,21},\quad E_3^{7,21} \cong E_2^{7,21}.$$
By assumption, the value of $Sq^1_3$ agrees with $Sq^1_2$ under the correspondences established by the last two of these isomorphisms.
So the diagram below appears on the $E_3$-page.
\begin{center}
\begin{tikzpicture}[scale=10, node distance=3cm]
\node (A) {$(6,21)\,\, \ecirc, \xcirc$};
\node (B) [right=3cm of A] {$(7,21)\,\, \xcirc$};
\node (C) [below of=A] {$(3,17)\,\, \xcirc$};
\node (D) [below of=B] {$(4,17)\,\, 0$};
\draw[->,double, bend left=20] (A) to node {$Sq^1_3\ne 0$} (B);
\draw[->,double, bend right=20] (A) to node [swap] {$Sq^1_3\ne 0$} (B);
\draw[->,double] (C) to node {$Sq^1_3=0$} (D);
\draw[->,dashed] (C) to node {$d_3 \ne 0$} (A);
\draw[->,dashed] (D) to node [swap] {$d_3 =  0$} (B);
\end{tikzpicture} 
\end{center}
Again, $E_3^{4,17} =0$ implies that $d_3 Sq^1_3 = 0$. On the other hand, we
shall see that $Sq^1_3 d_3 \ne 0$. To understand this composition, first
observe that the lemma above implies that $d_3(\xcirc) = d_3 X_3(\ecirc) =
X_3 d_3(\ecirc) = \xcirc$, where $X_3 = H(X_2,d_2)$ acts non-trivially by
Eqn. \eqref{reducedeq}. Second, $Sq^1_3 = [Sq^1_2] = [Sq^1]$ is non-zero on
this generator by virtue of the computer computation. So $d_3 Sq^1_3 \neq
Sq^1_3 d_3$, which contradicts the assumption that $Sq^1$ gives rise to an
operation.
\end{proof}

The Cartan formula allows us to construct examples for
which $Sq^n$ is not an operation on the Szab\'{o} spectral sequence.

\begin{corollary}
For each $n\geq 1$, there is link $K_n$ for which the Steenrod operation $Sq^n$ does not gives rise to an operation on Szab\'{o}'s geometric spectral sequence.
\end{corollary}
\begin{proof}
Set $K_n := \coprod_{i=1}^n T(4,5)$ so that $CKh(K_n) \cong CKh(T(4,5))^{\otimes n}$. The Szab\'{o} differential $\d_{Sz}$ respects this isomorphism and the K\"{u}nneth formula shows $E_m(K_n) \cong E_m(T(4,5))^{\otimes n}$ and
$E_{m+1}(T(4,5))^{\otimes n} \cong H(E_m(T(4,5))^{\otimes n}, d_m)$. 

Now let $a\in E^{3,17}_3$ be the element which satisfies $d_3 a=b \in E^{6,21}_3$ as in the proof of the previous theorem. In what follows, set $Sq^n := Sq^n_3$. We have $Sq^n a =0$ for all $n\geq 1$ and $Sq^nb = 0$ for all $n>1$. Consider that 
$d_3(a \ott b \ott \cdots \ott b) = b \ott \cdots \ott b$ and the Steenrod operation is
$$Sq^n(b\ott \cdots \ott b) = \sum_{i_1 + \cdots + i_n = n} Sq^{i_1} b \ott \cdots \ott Sq^{i_n} b = Sq^1 b \ott \cdots \ott Sq^1 b$$
which ensures that the composition $Sq^n d_3(a\ott b\ott \cdots \ott b)$ is non-zero. On the other hand, $Sq^n(a \ott b\ott \cdots \ott b) = Sq^1 a \ott Sq^1 b \ott \cdots \ott Sq^1 b = 0$ because $Sq^1a = 0$. Therefore, $d_3 Sq^n \ne Sq^n d_3$.
  \end{proof}

\section{Szab\'{o} homology as incompatible with Khovanov spectrification}
This section describes a setting in which the main theorem constitutes an
obstruction to a spectrification construction. Briefly, a spectrification
$\aY$ of Szab\'{o} homology cannot admit a filtration from which the
Lipshitz-Sarkar spectrification $\aX$ of (even) Khovanov homology can be
recovered in such a way as to be compatible with the spectral sequence
associated to a filtered stable homotopy type. This is because such a
filtration implies that the Steenrod operations on $H^*(\aX)$ extend
to the associated spectral sequence. 

As the literature in this area is not yet mature, the content of this
section should be read with caution.

R. Lipshitz and S. Sarkar's spectrification of Khovanov homology \cite{LS} can be described \cite{LLS} as follows: for each link diagram $K$, there is a wedge product of finite CW spectra
$$\aX(K) := \bigvee_{q\in\ZZ} \aX^q(K)$$
such that
\begin{enumerate}
\item For each $q\in \ZZ$, the cellular chain complex $C^*_{cell}(\aX^q(K))$
  is isomorphic to the Khovanov chain complex $CKh^{*,q}(K)$ with $q$-degree
  $q$. The isomorphism takes cells to monomials.
\item For each $q\in \ZZ$, the stable homotopy type of the spectrum
  $\aX^q(K)$ is an invariant of the link represented by $K$.
  \end{enumerate}

The definitions below are motivated by this characterization together with
the desire to obtain a reasonable looking obstruction statement. (Note that
we use cospectra $\aX, \aY\in Ob(Sp^{op}) = Ob(Sp)$
because cohomology is a contravariant homological functor.)

\begin{definition}\label{compatdef}
Suppose that $(C,d)$ is a chain complex of abelian groups, a {\em spectrification} $(C,d)$ is a cospectrum $\aX\in Sp^{\op}$ for which there is a canonical isomorphism: 
$$H^*(\aX;\ZZ) \xto{\sim} H(C,d)$$
between the ordinary cohomology $H^*(\aX;\ZZ) = Hom_{Sp^\op}(H\ZZ,\aX)$ of $\aX$ and the homology of the chain complex $(C,d)$. (In the knot theory context, a map is {\em canonical} when it commutes with maps induced by homotopy equivalences among spectra associated to Reidemeister moves.)

Suppose $\aX$ is a spectrification of $(C,d)$. Given a perturbation $\d : C \to C$, and a filtration 
$$\cdots \subseteq F^nC \subseteq F^{n-1} C \subseteq \cdots $$
with identifications $\bigcup_{n\in\ZZ} F^{n} C \cong (C,d+\d)$ and
$$\bigoplus_{n\in \ZZ} F^{n}C/F^{n-1} C \xto{\sim} (C,d).$$
A spectrification $\aY \in Sp^{\op}$ of a chain complex $(C, d + \d)$ is {\em an extension of a spectrification $\aX$ of the chain complex $(C, d)$ along the perturbation $\d$} when $\aY\in Sp^{\op}$ is a filtered cospectrum with associated graded $\aX$: in more detail, there is a functor $F : (\ZZ,<)\to Sp^{\op}$ of $(\infty,1)$-categories corresponding to assignments
$$\cdots \to F^{n-1} \aY \xto{\kappa_n} F^{n} \aY \to \cdots$$
with canonical equivalences
$\colim F^n \aY\cong \aY$  and 
$$E_1^{p,q}= H^{p+q}(\cofib(\kappa_p)) \xto{\sim} C^p_{cell}(\aX^q).$$
\end{definition}

\begin{remark}
A filtered cospectrum $\aY\in Sp^{\op}$ and a homology functor $H : Sp^{\op}\to Ab$, such as $H = Hom_{Sp^{\op}}(H\FF_2, \cdot)$, give rise to a spectral sequence. The construction of this spectral sequence is a lift of the spectral sequence associated to a filtered chain complex. Our definitions above have been chosen to be compatible with the reference \cite[\S 1.2.2]{Lurie}. 
\end{remark}

The proposition below observes that cohomology operations commute with this construction in our setting.

\begin{prop}
Suppose $\aX$ is a spectrification of $(C,d)$. Then a stable cohomology operation $Sq^n : H^i(\aX) \to H^{i+n}(\aX)$ gives rise to an operation on the spectral sequence $\{ E_r \}_{r\geq 2}$ associated to any extension $\aY$ of this spectrification.
\end{prop}
\begin{proof}
The $E_r$-page is defined in Construction \S 1.2.2.6 \cite{Lurie} to be the image
  \begin{equation*}
    E_r^{p,q} := im(H^{p+q}(\aY(p-r,p)) \xto{(\kappa^{p,q}_r)^*} H^{p+q} (\aY(p-1, p+r-1)))
\end{equation*}    
of a map induced on cohomology between auxilliary cospectra of the form $\aY(p,p+m) = \cofib(F^p\aY\to F^{p+m}\aY)$ in $Sp^{\op}$. Naturality of $Sq^n$, $Sq^n(\kappa^{p,q}_r)^* = (\kappa^{p,q}_r)^* Sq^n$, implies that $Sq^n : im(\kappa^{p,q}_r)^* \to im(\kappa^{p,q}_r)^*$ since if $x\in im(\kappa^{p,q}_r)^*$ then $Sq^n x = Sq^n (\kappa^{p,q}_r)^* y = (\kappa^{p,q}_r)^* (Sq^n y)$. So $Sq^n$ defines a map on the $E_r$-page.  

In a similar way, the map $Sq^n$ commutes with the differential $d_r : E_r^{p,q} \to E_r^{p-r,q+r-1}$ because $Sq^n$ is a stable operation and the differential is uniquely determined a commutative diagram involving connecting homomorphisms \cite[p. 49]{Lurie}.

The map $Sq^n$ on $E_r$ is induced by the map $Sq^n$ on $E_{r-1}$ by construction.
  \end{proof}

\begin{corollary}
There is no extension of the (even) spectrification $\aX$ of Khovanov
homology to a spectrification $\aY$ of Szab\'{o} homology.
  \end{corollary}
\begin{proof}
The main theorem above shows that the Steenrod operations on $H^*(\aX)$ determined by the spectrification of (even) Khovanov homology do not extend to operations on the associated Szab\'{o} spectral sequence. This is inconsistent with the proposition above.
  \end{proof}

\begin{remark}
The observations above do not obstruct an approach to the spectrification of Szab\'{o} homology extending the recent spectrification of odd Khovanov homology \cite{SSS}.
  \end{remark}

\section{Integral lifts}
The main theorem in this section shows that there is no $\ZZ$-lift of the
$\FF_2$-Szab\'{o} spectral sequence from the (even) Khovanov homology. Such
a lift is inconsistent with a Bockstein operation on Khovanov homology which
extends to an operation on the Szab\'{o} spectral sequence. This is the
content of the technical lemma below. However, before we reach the lemma,
we'll need some definitions.

\begin{definition}
  If $(C,d)$ is a $\FF_2$-chain complex then a  \emph{$\ZZ$-lift $(\tilde{C},\tilde{d})$ of $(C,d)$} is a chain complex of free abelian groups and an isomorphism of chain complexes
  $$\varphi : (\tilde{C}, \tilde{d}) \otimes_{\ZZ} \FF_2 \xto{\sim} (C,d).$$
\end{definition}

\begin{lemma}
If $(\tilde{C}, \tilde{d})$ is $\ZZ$-lift of $(C,d)$ then there is a surjective {\em $\FF_2$-reduction map} $\rd{x} := \varphi(x\ott 1)$ and the (not unique) {\em $\ZZ$-lift map} $x\mapsto \tilde{x}$ 
$$\rd{\cdot} : (\tilde{C},\tilde{d}) \rightleftarrows (C,d) : \tilde{\cdot}$$
which satisfy
\begin{enumerate}
\item $\rd{\tilde{x}} = x$, 
\item $d\rd{x} = \rd{\tl{d}x}$ and
\item $[\rd{x}] = \rd{[x]}$
\end{enumerate}
for all $x\in C$. 
  \end{lemma}
\begin{proof}
The first by definition of lift, the second because $\tl{d}\equiv d\,\, (\mathrm{mod}\,\, 2)$ and, for the last one, $\FF_2$-reduction commutes with passage to homology because the second implies that $\FF_2$-reduction is a chain map.
\end{proof}
 
The Bockstein $\beta : H^n(C,d) \to H^{n+1}(C,d)$ associated to an
$\FF_2$-chain complex $(C,d)$ agrees with the Steenrod operation $Sq^1$, \cite[p. 489 (7)]{Hatcher}. It has a relatively simple formula, in the presence of a $\ZZ$-lift $(\tilde{C},\tilde{d})$, as the $\FF_2$-reduction \cite[Ex. 3.2]{Labbockstein} of the connecting map for the long exact sequence associated to the short exact sequence
$$0\to \tl{C}\xto{2\cdot} \tl{C} \to C \to 0.$$

\begin{definition}
  If $(\tilde{C}, \tilde{d})$ is $\ZZ$-lift of $(C,d)$ and $x\in C$ is a $d$-cycle then the Bockstein is given by $\beta(x) = \rd{\frac{1}{2} \tilde{d} \tilde{x}}$. 
  \end{definition}

The definitions below will be used by the Technical Lemma which follows.

\begin{definition}\label{compdef}
  Suppose $(\tilde{C},\tilde{d})$ is a $\ZZ$-lift of $(C,d)$ and $\{ F^nC \}_{n\in\ZZ}$
is a filtration of $(C,d)$. Then a filtration $\{F^n \tilde{C}\}_{n\in\ZZ}$ of the $\ZZ$-lift $(\tilde{C},\tilde{d})$ is {\em compatible} with the filtration $\{F^n C\}_{n\in \ZZ}$ when the two properties below are satisfied.
\begin{enumerate}
\item $F^n\tilde{C}$ is a $\ZZ$-lift of $F^nC$, i.e. there are maps
  $$\varphi_n : (F^n\tilde{C},\tilde{d})\otimes_{\ZZ} \FF_2 \xto{\sim} (F^nC, d)$$
\item   The maps $\varphi_n$ preserve the filtrations
$$\begin{tikzpicture}[scale=10, node distance=1.75cm]
\node (A) {$F^n\tilde{C}\ott \FF_2 $};
\node (C) [right=.6cm of A] {$F^{n-1}\tilde{C}\ott \FF_2$};

\node (B) [below of=A] {$F^nC$};

\node (D) [below of=C] {$F^{n-1}\tilde{C}$};
\draw[->] (A) to node {$\subseteq$} (C);
\draw[->] (B) to node {$\subseteq$} (D);

\draw[->] (A) to node [swap] {$\varphi_n$} (B);
\draw[->] (C) to node {$\varphi_{n-1}$} (D);
\end{tikzpicture}$$  
  \end{enumerate}
Such a compatible filtration $\{ F^n\tilde{C} \}_{n\in\ZZ}$ is {\em torsion-free} when the associated graded $\oplus_{n\in\ZZ} F^n \tilde{C}/F^{n-1}\tilde{C}$ is torsion-free.
\end{definition}

It is time to state the technical lemma.

\newcommand{\hrt}{\ensuremath{(\bigstar)}}
\begin{lemma}{(The Technical Lemma)}
Suppose that $(C, d)$ is a chain complex of $q$-graded $\FF_2$-vector spaces, $C = \oplus_{(t,q)\in \ZZ\times\ZZ} C^{t,q}$ and  $d = \sum_{i=1}^\infty d_i$, where $d_i$ of $(t,q)$-bidegree $(i,1)$. The filtration given by $F^pC:=\{ x : \vnp{x}_t \geq p \}$  so that
$$\cdots \subseteq F^p C \subseteq F^{p-1}C \subseteq \cdots$$
determines a spectral sequence $\{E_r,\d_r\}_{r=2}^{\infty}$ with $E_2^{p,q} \cong H(C, d_1)$. If $(\tilde{C}, \tilde{d})$ is a $\ZZ$-lift of $(C, d)$ with a compatible and torsion-free filtration, then the Bockstein of $(C,d_1)$ extends to an operation $\{\beta_r : E^{p,q}_r \to E^{p+1,q}_r\}_{r \geq 2}$ on the spectral sequence. 

Moreover, the maps $\beta_r$ are determined by their values on $[x] \in E^{p,q}_r$ by choosing a representative $x = \sum_i x_i \in [x]$ with $(t,q)$-bidegree  $\vnp{x_i}_{(t,q)}=(p+i,p+q)$ and setting
$$  \beta_r([x]) :=     [\rd{\frac{1}{2} \tilde{d}_1 \tilde{x}_0} ].$$
 \end{lemma}
 \begin{proof}
   Recall from the construction of the spectral sequence, as in the proof of Thm. 2.6 \cite{McCleary}, that the $E_r$-page consists of vector spaces of the form
   $$E^{p,q}_r = Z_r^{p,q}/(Z_{r-1}^{p+1,q-1} + B_{r-1}^{p,q})$$
   where the almost cycles and almost boundaries are given respectively by
$$Z_r^{p,q} := F^pC^{p+q}\cap d^{-1}(F^{p+r} C^{p+q+1})\conj{ and } B_r^{p,q} := F^pC^{p+q}\cap d(F^{p-r} C^{p+q-1}).$$
By assumption the structure of this spectral sequence is compatible with the choice of integral lift.
Properties (1) and (2) from Def. \ref{compdef} imply that 
 there are operations $$\rd{\cdot} : \{\tl{E}^{p,q}_r, \tl{\d}_r\}_{r=2}^\infty \rightleftarrows \{E^{p,q}_r, \d_r\}_{r=2}^\infty : \tl{\cdot}.$$

   Now in order to introduce the Bockstein spectral sequence operation, observe that, in our setting, almost cycles can be written more concretely as
   $$Z_r^{p,q} = \{ \sum_{i\geq 0} x_i : \substack{\vnp{x_i}_t = p+i\\\vnp{x_i}_q = p+q} \conj{ and } \hrt \sum_{\substack{i+j=N\\j\geq 0,i\geq 1}} d_i(x_j) = 0 \normaltext{ for } 1\leq N < r  \}.$$
set  $P^{p,q}_r := Z_r^{p,q}/Z_{r-1}^{p+1,q-1}$. Since $r\geq 1$, $\hrt$ implies that $x=\sum_i x_i\in Z_r^{p,q}$ satisfies $d_1(x_0) = 0$, the definition
$$\beta_r : P^{p,q}_r \to F^{p+1}C^{p+1+q} \conj{ given by }   \beta_r(x_0) := \rd{\frac{1}{2} \tilde{d}_1 \tilde{x}_0}$$
is in agreement with the $d_1$-Bockstein. 

The claims (i) and (ii) below show that $\beta_r$ is well-defined and descends to a map $\beta_r : P^{p,q}_r \to P^{p+1,q}_r$ respectively. Claim (iii) shows that $\beta_r$ descends further to a map $\beta_r : E^{p,q}_r \to E^{p+1,q}_r$. Finally, claims (iv) and (v) show that the maps $\beta_r$ determine an operation on the spectral sequence.
\begin{enumerate}[label=(\roman*)]
\item $\beta_r(Z_{r-1}^{p+1,q-1})=0$: We claim that if $x, x'\in Z^{p,q}_r$ are of the form $x = x_0 + \sum_{i\geq 1} x_i$ and $x' = x_0 + \sum_{i\geq 1} x'_i$ where $\vnp{x_i}_{(t,q)} = (p+i, p+q)$ and $\vnp{x'_i}_{(t,q)} = (p+i, p+q)$ then $x - x' \in Z^{p+1,q-1}_{r-1}$. This follows from $\hrt$: 
$$(dx)_n=\sum_{i+j=n} d_i x_j = 0 \conj{ and } (dx')_n=\sum_{i+j=n} d_i x'_j = 0$$
imply $(d(x-x'))_n=\sum_{i+j=n} d_i (x_j-x'_j) = 0$, so $ x-x'=\sum_{i\geq 1} (x_i-x'_i) \in Z^{p+1,q-1}_r$.
\item $im(\beta_r) \subseteq Z^{p+1,q}_r$: Suppose $x = x_0 + x_1 + \cdots$, choose lifts $\tilde{x} = \tilde{x}_0 + \tilde{x}_1 +\cdots$ the relation $\hrt$ becomes
  $$(\tilde{d}\tilde{x})_n = \sum_{i+j=n} \tilde{d}_i \tilde{x}_j =: 2 \tilde{y}_{n-1}.$$
  Set $y_n := \rd{\tilde{y}_n}$. Since $\beta_r(x) = y_0$, in order to prove the claim, it suffices to show that $y_0 + y_1 + \cdots$ satisfies $\hrt$. This is implied by the corresponding relation for the $\ZZ$-lift:
  $$(\tilde{d}\tilde{y})_n = \sum_{i+j=n} \tilde{d}_i \tilde{y}_j = \sum_{m=0}^{n-1} \frac{1}{2} \left( \sum_{i+j=n-m} \tilde{d}_i \tilde{d}_j\right) \tilde{x}_m = 0$$
which is zero because $\tilde{d}^2 = 0$ implies $\sum_{i+j=n-m} \tilde{d}_i \tilde{d}_j = 0$.
\item $\beta_r(B_{r-1}^{p,q})=0$: If $dy = x$ then $\tilde{x}_0=(\tilde{d}\tilde{y})_0$ which implies that $\tilde{d}_1 \tilde{x}_0 = \tilde{d}_1(\tilde{d}\tilde{y})_0$. Now  $(\tilde{d}(\tilde{d}\tilde{y}))_1 = 0$ because $\tilde{d}^2=0$, on the other hand
$$(\tilde{d}(\tilde{d}\tilde{y}))_1 =  \sum_{i+j=1} \tilde{d}_i (\tilde{d}\tilde{y})_j= \tilde{d}_1 (\tilde{d}\tilde{y})_0 + 0 + \cdots + 0$$
the first equality is the degree $1$ term of $\tilde{d}(\tilde{d}\tilde{y})$ and the second equality holds because $(\tilde{d}_n(\tilde{d}\tilde{y})_{-n+1})= 0$ for $n \neq 0$ because $\tilde{d}_n = 0$ for $n \leq 0$ and $(\tilde{d}\tilde{y})_{-n+1} = 0$ for $n > 0$. So $\tilde{d}_1 (\tilde{d}\tilde{y})_0 = (\tl{d}(\tl{d}\tl{y}))_1 =  0$. All together now, 
$$\beta_r(dy) = \rd{\frac{1}{2} \tilde{d}_1 (\tilde{d}\tilde{y})_0} = \rd{\frac{1}{2}\cdot 0} = 0.$$

\item $\beta_r \d_r = \d_r \beta_r$: From \cite[p. 35]{McCleary}, a formula for the differential 
  $$\d_r([x]) = [(dx)_r] = [ \sum_{i+j=r} d_i x_j] \conj{ for } x = \sum_i x_i \in Z^{p,q}_r$$
of the spectral sequence can be observed from the definition of $\d_r$ as $d=\sum_i d_i$ in the image of the quotient $\eta : Z^{p,q}_r \to E^{p,q}_r$.

Now when $(dx)_n =0$, there are elements $\tilde{y} = \tilde{y}_0 + \tilde{y}_1 + \cdots$ so $(\tilde{d}\tilde{x})_n = 2\tilde{y}_{n-1}$  so that 
  $$\beta_r([x_0]) = [\rd{\frac{1}{2} \tilde{d}_1 \tilde{x}_0}] =[\rd{\frac{1}{2} (\tilde{d}\tilde{x})_1} ] = [\rd{\tilde{y}_0}]=[y_0].$$
  Now on one hand $(\tilde{d}(\tilde{d}\tilde{x}))_{r+1} = 2(\tilde{d}\tilde{y})_r$, but on the other hand $(\tilde{d}(\tilde{d}\tilde{x}))_{r+1}$ is
        \setlength{\multlinegap}{50pt}
\begin{multline*}
  \tilde{d}_1(\tilde{d}\tilde{x})_{r} + \tilde{d}_2(\tilde{d}\tilde{x})_{r-1} + \cdots + \tilde{d}_r(\tilde{d}\tilde{x})_{1}\\
  = \tilde{d}_1(\tilde{d}\tilde{x})_{r} + 2\tl{d}_2(\tl{y}_{r-2}) + 2\tl{d}_3(\tl{y}_{r-3}) + \cdots + 2\tl{d}_r(\tl{y}_0).
  \end{multline*}
It follows that $\rd{\tilde{d}_1(\tilde{d}\tilde{x})_r} = \rd{\tl{d}(\tl{d}\tl{x})_{r+1}} = \rd{2 (\tilde{d}\tilde{y})_r}$.
So
 \begin{multline*}
   (\beta_r \d_r)([x_0]) = \beta_r([(dx)_r])
   = [\rd{\frac{1}{2} \tilde{d}_1 (\tilde{d}\tilde{x})_r}]= [\rd{\frac{2}{2} (\tilde{d}\tilde{y})_r}]
   =\rd{[(\tilde{d}\tilde{y})_r]}\\
   = \rd{\tl{\d}_r [\tilde{y}_0]}= \d_r \rd{[\tl{y}_0]}
   = \d_r [\rd{\tl{y}_0}]
   = \d_r [y_0]
   = (\d_r \beta_r)([x_0]).
   \end{multline*}
\item $\beta_{r+1} = H(\beta_r, d_r)$: By construction, $\beta_{r+1}([x]) = [\beta_{r}(x)]$.
\end{enumerate}

   \end{proof}

The above lemma applies in the Szab\'{o} setting by taking $C := CKh(K;\FF_2)$, $d := \d_{Sz}$ the Szab\'{o} differential and $d_1 := d_{Kh}$ the Khovanov differential. The differential $\d_{Sz}$ is degree $+1$ in the $\delta$-grading, $\delta(x) := \vnp{x}_t - \vnp{x}_q/2$. In agreement with Szab\'{o} \cite[Def. 7.1]{Szabo}, the filtration is given by
$$F^p C^q := \{ x \in C : \vnp{x}_\d = q, \vnp{x}_t \geq p \}.$$

 \begin{theorem}
There is no integral lift of the Szab\'{o} spectral sequence for which the $E_2$-page agrees with (even) Khovanov homology. 
  \end{theorem}
 \begin{proof}
Suppose there is such a chain complex $(CKh(T(4,5);\ZZ),\tilde{\d}_{Sz})$, so
$$  (CKh(T(4,5);\FF_2), \d_{Sz}) \cong (CKh(T(4,5);\ZZ),\tilde{\d}_{Sz}) \ott_{\ZZ} \FF_2$$
and $\tilde{\d}_{Sz} = d_1 + d_2 + \cdots$ where $d_1$ is the usual (even) Khovanov differential. By setting $F^p CKh(T(4,5);\ZZ) := \{ x : \vnp{x}_\d = q, \vnp{x}_t \geq p \}$, we obtain a filtration which is compatible and torsion-free. Then, by the preceeding lemma, the Bockstein $\beta$ associated to $\tilde{\d}_{Sz}$ must give an operation on the Szab\'{o} spectral sequence which agrees with the Bockstein $\beta$ associated to the (even) Khovanov $d$ differential on the $E_2$-page. The main theorem above shows that this is not possible.
  \end{proof}

\begin{remark}
In contrast to the theorem above there is a preprint suggesting that
Szab\'{o}'s spectral sequence lifts to {\em odd} integral Khovanov homology
\cite{Beier}.  If this is true then then it follows from the technical lemma that 
the odd Bockstein map extends to an operation on the Szab\'{o} spectral sequence. The even Bockstein and odd Bockstein are known to be different, see \cite[Ex. 3.1]{PS}.
\end{remark}


\section{A pattern}
In the first few dozen knot table examples one can observe a relationship
between $Sq^2$ and $d_2$. It appears that every non-zero $d_2$ leads to a non-zero $Sq^2$,
they occur together in the butterfly configuration depicted below.

\begin{center}
\begin{tikzpicture}[scale=1]
  \draw[->] (0,0) -- (3.5,0) node[right] {$t$};
  \draw[->] (0,0) -- (0,3.5) node[above] {$q$};
  \draw[step=1] (0,0) grid (3,3);
  \draw (0.5,-.2) node[below] {$t$};
  \draw (1.5,-.2) node[below] {$t+1$};
  \draw (2.5,-.2) node[below] {$t+2$};
  \draw (3.5,-.2) node[below] {};
  \draw (-.2,0.5) node[left] {$q$};
  \draw (-.2,1.5) node[left] {$q+2$};
  \draw (-.2,2.5) node[left] {$q+4$};
  \draw (-.2,3.5) node[left] {};

\fill (0.5, 0.5) circle (.15);
  \draw (0.5, 1.5) circle (.15);

  \fill (2.5, 1.5) circle (.15);
  \draw (2.5, 2.5) circle (.15);
\draw[->,double] (0.750, 1.5) -- (2.25, 1.5);
  \draw[->] (0.768, 0.634) -- (2.232, 1.366);
  \draw[->] (0.768, 1.634) -- (2.232, 2.366);

\end{tikzpicture}
\end{center}
See for example, $(t,q) = (2, 15)$ in $T(4,5)$.  More formally, this pattern suggests that
$$d_2(xm) = xm' \conj{ implies } Sq^2(m) = xm' + \cdots.$$ 
The converse is false as $T(4,5)$ contains a $Sq^2$-operation at $(t,q) =(5,21)$ of a different sort. 


\subsection*{Acknowledgments} The authors thank R. Lipshitz for mentioning this
question at the Banff workshop in June 2021, the anonymous referee and
E. Magnuson for her illustrations.  The first author thanks the organizers
M. Aganagi\'{c}, S. Krushkal and B. Webster for inviting his
participation. This paper was partially funded by Simons Award \#638089.

\bibliography{spectralnote}{}
\bibliographystyle{amsalpha}

\end{document}